\documentclass[11pt]{article}
\usepackage{latexsym, amsfonts, amsmath, amssymb, amscd, epsfig}
\usepackage[numbers,sort&compress]{natbib}
\usepackage{enumerate}
\usepackage{amscd}
\usepackage{array}
\usepackage{amsfonts}
\usepackage{amsthm}
\usepackage{verbatim}
\usepackage{mathrsfs}
\usepackage{enumerate}
\usepackage{color}

\numberwithin{equation}{section}
\newtheorem{thm}{Theorem}[section]

\newtheorem{lem}[thm]{Lemma}

\newtheorem{cor}[thm]{Corollary}

\newenvironment{pf}{{\noindent \it \bf Proof:}}{{\hfill$\Box$}\\}

\newcommand{\be}{\begin{equation}}
\newcommand{\ee}{\end{equation}}
 \newcommand{\ba}{\begin{aligned}}
 \newcommand{\ea}{\end{aligned}}

\newcommand{\R}{\mathbb R}

\newcommand{\B}{\Big}
  \newcommand{\f}{\frac}
\newcommand{\norm}[1]{|\!|#1|\!|}

\begin{document}

\title{Regularity criteria of 3D generalized magneto-micropolar fluid system in terms of the pressure}
\author{Jae-Myoung Kim$^{1}$\\
$^{{\small 1}}${\small Department of Mathematics Education, Andong National University,}\\
{\small  Andong 36729, Republic of Korea,}
\\ \small Email: jmkim02@anu.ac.kr}%
\date{}
\maketitle

\begin{abstract}
This work focuses on regularity criteria of 3D generalized magneto-micropolar fluid system in terms of the pressure in Lorentz spaces. 
\end{abstract}

\noindent {\textit{Mathematics Subject Classification(2000):}\thinspace
\thinspace 35Q35,\thinspace \thinspace 35B65,\thinspace \thinspace 76D05}
\newline
Key words: Regularity criteria; Weak solution, generalized magneto-micropolar fluid system \newline


\section{Introduction}

\bigskip This paper is concerned about the reglarity conditions of the weak
solutions to generalized magneto-micropolar fluid equations in
$\R^3$, which are described by
\begin{equation}
\label{mag-mipolar} \left\{
\begin{aligned}
\partial_t u +( u \cdot \nabla) u + \nabla \pi& =  -( \mu + \chi)\Lambda^{2\alpha} u+ \chi \nabla \times w + (b \cdot \nabla)  b,  \\
\partial_t w + (u \cdot \nabla) w & =  -\kappa \Lambda^{2\gamma} w+  \eta\nabla (\nabla \cdot w) + \chi \nabla \times u -2 \chi w,  \\
\partial_t b + (u \cdot \nabla) b & =  -\nu \Lambda^{2\beta} b+ (b \cdot \nabla) u, \\
\nabla \cdot u (\cdot,t) & =  \nabla \cdot b(\cdot,t)\, = 0,
\end{aligned}
\right.
\end{equation}
where $u = u(x, t)$, $w = w(x, t)$, $b = b(x, t)$ and $\pi: =\pi(x,
t)=\mathcal{P}+\frac{|b|^2}{2}$ denote the fluid velocity, the micro-rotation velocity (angular
velocity of the rotation of the fluid particles), the magnetic and
total pressure fields respectively. The notation $ \Lambda :=
(-\Delta)^{1/2}$ stands for positive constants. The positive constant $\kappa$
in (\ref{mag-mipolar}) correspond to the angular viscosity, $\nu$ is
the inverse of the magnetic Reynolds number and $\chi$ is the
micro-rotational viscosity. We consider the initial value problem of
\eqref{mag-mipolar}, which requires initial conditions
\begin{equation}\label{ini}
u(x,0)=u_0(x),\quad w(x,0)=w_0(x) \quad \text{and} \quad
b(x,0)=b_0(x), \qquad x\in\R^3,
\end{equation} and we also assume that $ \text{div} \
u_0=0=\text{div} \ b_0$. 

The authors in \cite[Theorem 2.2]{CH19} construct the
existence of Leray-Hopf solution on the whole space $\R^3 \times (0, T)$ to the generalized Navier-Stokes equation. It is worth pointing that
$\dot{H}^{\frac{5-4\alpha}{2}}$ is a critical space, that is,
$\dot{H}^{\frac{5-4\alpha}{2}}$ norm is scaling invariant  is also a solution, where
$u_\lambda(x, t)=\lambda^{2\alpha-1}u(\lambda x,
\lambda^{2\alpha}t)$ and $\pi_\lambda(x, t)=\lambda^{2(2\alpha-1)}\pi(\lambda x,
\lambda^{2\alpha}t)$ with any $\lambda>0$ . By
Sobolev embedding theorem, it is checked that
$\dot{H}^{\frac{5-3\alpha}{2}}(\R^3)\hookrightarrow
L^{\frac{6}{3-2\alpha}}(\R^3)$. Then
\[
\|u\|_{L^4(0,T ;L^{\frac{6}{3-2\alpha}})} \leq
C\|u\|^2_{L^\infty(0,T ;\dot{H}^{\frac{5-4\alpha}{2}})}
\|u\|^2_{L^2(0,T ;\dot{H}^{\frac{5-2\alpha}{2}})},\quad
\frac{2\alpha}{4}+\frac{3}{\frac{6}{3\alpha-2}}=2\alpha-1.
\]


For the regularity issues to weak solutions to the generalized Navier-Stokes equation, we refer to \cite{Chae07}, \cite{FO08}, \cite{FFY13}.

Recently, Deng and Shang \cite{DS21} obtained global-in-time existenceand uniqueness of smooth solutions to the problem \eqref{mag-mipolar}--\eqref{ini} if $\alpha\geq \frac{1}{2}+\frac{n}{4}, \alpha +\gamma \geq \max (2,\frac{n}{2}), \mbox{and}\quad \alpha +\beta \geq 1 +\frac{n}{2}$. On the other hands, Fan and Zhong \cite{FS22} established local-in-time existene and uniqueness of smooth solutions to the problem \eqref{mag-mipolar}--\eqref{ini} for  $\alpha+\gamma>1$ and furthermore they gave some regularity criteria via  the gradient of velocity in  a meaningful and appropriate spaces.

For the regularity criteria in Lorentz space, Li
and Niu \cite{LN22} proved that a weak solution $(u,\,w,\ b)$ for the standard 3D
MHD equations become regular under the scaling invariant conditions for the total pressure, in particualr,
so called Serrin's conditions,
$\pi \in L^{q,\infty}(0,\, T; L^{p,\infty}({\mathbb{R}}^3))$
with $3/p+2/q \leq 2$ and $p>\frac{3}{2}$ (compare to 
\cite{[BG]}, \cite{Zhou05}, \cite{Zhou06}, \cite{Zhou07}, \cite{BPR14}, \cite{[Suzuki1]}, \cite{[Suzuki2]} for Navier--Stokes equations). 
%
%
%
For
$p,q\in[1,\infty]$, we define
$$
\|f\|_{L^{p,q}(\Omega)}=\left\{\ba
&\B(p\int_{0}^{\infty}\alpha^{q}|\{x\in \Omega:|f(x)|>\alpha\}|^{\frac{q}{p}}\frac{d\alpha}{\alpha}\B)^{\frac{1}{q}} , ~~~q<\infty, \\
 &\sup_{\alpha>0}\alpha\ |\{x\in \Omega:|f(x)|>\alpha\}|^{\frac{1}{p}} ,~~~q=\infty.
\ea\right.
$$
 (see e.g. \cite{Neil}, \cite{Grafakos}, \cite{Maly})
Motivated by the recent works in \cite{FS22} and \cite{LN22}, the purpose of this note is to
establish regularity criteria of 3D generalized magneto-micropolar fluid system \eqref{mag-mipolar} in terms of the pressure in Lorentz space. In this paper, we assume
that {\color{blue}$\mu+\chi=1$, $\nu=\eta=1$}  and $\alpha=\beta$.

%

Our results are stated as follows.

\begin{thm}\label{the1.1}
Let  $\quad 0 < T < \infty$ and $u_0, b_0, w_0   \in  H^m(\R^3)$ with $m > \frac{5}{2}$
and $1\leq \alpha, \gamma \leq \frac{5}{4}$. There
exists a sufficient constant $\epsilon
> 0$ such that if $\pi$ or $\nabla \pi$ satisfy
\begin{enumerate}[(A)]
 \item   $  \pi \in L^{q,\infty}(0,T; L ^{p,\infty}(\mathbb{R}^{3}))$ ~and~ $$\|\pi\|_{L^{q,\infty}(0,T; L ^{p,\infty}(\mathbb{R}^{3}))}
     \leq\varepsilon, ~ \text{with} ~~\frac{3}{p}+\frac{2\alpha}{q}=2(2\alpha-1) , ~\frac{3}{2(2\alpha-1)}<p<\infty;  $$
 \item  $ \nabla \pi \in L^{\frac{2r\alpha}{3r\alpha-3}}(0,T; L ^{r,\infty}(\mathbb{R}^{3}))$ with $  \frac{1}{\alpha} < r \leq\infty,
$
 \end{enumerate}
 then a weak $\left( u,b\right) $ is regular on $(0,T].$
\end{thm}

%

\begin{cor}
Let $\alpha=\beta=\gamma=1$. Assume $(u_0, b_0) \in L^2_{\sigma}(\mathbb{R}^{3}) \cap L^4_{\sigma}(\mathbb{R}^{3})$  and $w_0 \in  L^2(\mathbb{R}^{3}) \cap L^4(\mathbb{R}^{3})$. 
Let the triple $(u, b, w)$ be a weak solution
to system \eqref{mag-mipolar} on some time interval $[0, T)$ with $0<T<\infty$. There exists a sufficient constant $\epsilon
> 0$ such that
suppose that $\nabla \pi$ satisfies $ \nabla \pi \in L^{\frac{2r}{3r-3},\infty}(0,T; L ^{r,\infty}(\mathbb{R}^{3}))$  with
\[
\| \nabla \pi\|_{L^{\frac{2r}{3r-3},\infty}(0, T
;L^{r,\infty})}\leq \epsilon, \quad 1 < r \leq \infty,
\]
 then the weak solutions $\left( u,b, w\right) $ is regular on $(0,T].$
\end{cor}

\begin{thm}\label{the1.2}
Let $\alpha=\beta=\gamma=1$. Assume $(u_0, b_0) \in L^2_{\sigma}(\mathbb{R}^{3}) \cap L^4_{\sigma}(\mathbb{R}^{3})$  and $w_0 \in  L^2(\mathbb{R}^{3}) \cap L^4(\mathbb{R}^{3})$. 
Let the triple $(u, b, w)$ be a weak solution
to system \eqref{mag-mipolar} on some time interval $[0, T)$ with $0<T<\infty$.
There exists a sufficient constant $\epsilon
> 0$ such that if $\nabla \mathcal{\mathcal{P}}$  and $b$ satisfy $ \nabla \mathcal{\mathcal{P}} \in L^{\frac{2r}{3r-3},\infty}(0,T; L ^{r,\infty}(\mathbb{R}^{3}))$ and $ b \in L^{\frac{2a_1}{a_1-3},\infty}(0,T; L ^{a_1,\infty}(\mathbb{R}^{3}))$ with
\[
\| \nabla \mathcal{\mathcal{P}}\|_{L^{\frac{2r}{3r-3},\infty}(0, T
;L^{r,\infty})}\leq \epsilon, \quad 1 < r \leq \infty,
\]
and 
\[
\| b\|_{L^{\frac{2a_1}{a_1-3},\infty}(0, T
;L^{a_1,\infty})}< \infty, \quad 3 < a_1 \leq \infty,
\]
 then the weak solutions $\left( u,b, w\right) $ is regular on $(0,T].$
\end{thm}

\section{Proof of Theorem \protect\ref{the1.1}}

%
%

To control the fractional diffusion term, we recall the following
result (see e.g. \cite{Codorba}).

\begin{lem}\label{positivity}
With $0 < \alpha < 2$, $v,\Lambda^\alpha v \in L^p(\R^3)$ with
$p=2k$, $k\in \mathbb{N}$, we obtain
 \[
\int |v|^{p-2}v\Lambda^\alpha v\, dx \geq \frac{1}{p} \int
|\Lambda^{\frac{\alpha}{2}} v^{\frac{p}{2}}|^2\,dx.
 \]
\end{lem}

Also, we recall the following nonlinear Gronwall-type inequality
established in \cite{PX20} (see also \cite{BPR14} and \cite{LRM16}).
\begin{lem}\label{grow}
Let $T > 0$ and $\varphi \in L_{loc}([0, T ))$ be non-negative
function. Assume further that
\[
\varphi(t) \leq  C_0 + C_1\int^t_0 \mu(s)\varphi(s)\,ds + \kappa
\int^t_0 \lambda(s)^{1-\epsilon}\varphi(s)^{1+A(\epsilon)}\,
ds,\quad \forall \ 0 < \epsilon < \epsilon_0.
\]
Where $\kappa, \epsilon_0 > 0$ are constants, $\mu \in L^1(0, T )$
and $A(\epsilon)> 0$ satisfies $\lim_{\epsilon \rightarrow
0}\frac{A(\epsilon)}{\epsilon}= c_0 > 0$. Then $\varphi$ is bounded
on $[0, T ]$ if $\|\lambda\|_{L^{1,\infty}(0,T)} < c^{-1}_0
\kappa^{-1}$.
\end{lem}

In order to derive the regularity criteria of weak solutions to the system (\ref{mag-mipolar}), we introduce the definition of weak
solution. Let us denote
\begin{equation}\label{variable}
z^+=u+b,\ \ \text{\ }z^-=u-b.
\end{equation}%
Then system (1.1) can be reformulated as

\begin{equation}
\left\{
\begin{array}{c}
\partial _{t}z^{+}-\Lambda^{2\alpha}
z^{+}+(z^{-}\cdot \nabla )z^{+}-\chi\nabla \times w+\nabla \pi =0, \\
\partial _{t}z^{-}-\Lambda^{2\alpha}
z^{-}+(z^{+}\cdot \nabla )z^{-}-\chi\nabla \times w+\nabla \pi =0, \\
\partial _{t}w-\kappa\Lambda^{2\gamma}
w+ (\frac{z^++ z^-}{2})\cdot \nabla w - \nabla \mbox{div}\ w +
w-\nabla \times(\frac{z^++ z^-}{2}) =0, \\
\nabla \cdot z^{+}=\nabla \cdot z^{-}=0, \\
z^{+}(x,0)=z^{+}_{0}(x),\text{ \ }z^{-}(x,0)=z^{-}_{0}(x),%
\end{array}%
\right.  \label{fmmp-100}
\end{equation}

It is easy to show the following global $L^2$-bound,
\[
\|(u,b,w)((\tau))\|^2_{L^2}+ \int^t_0 (\|\Lambda^\alpha u(\tau )\|^2_{L^2} +
\|\Lambda^\gamma w(\tau)\|^2_{L^2}+ \|\Lambda^\alpha
b(\tau)\|^2_{L^2})\,d\tau \leq C.
\]

\begin{pf}
Multiplying the first and the second equations of \eqref{fmmp-100} by $%
\left\vert z^{+}\right\vert ^{2}z^{+}$ and $\left\vert
z^{-}\right\vert ^{2}z^{-}$, respectively, integrating by parts and
summing up, we have
\[
\frac{1}{4}\frac{d}{dt}\norm{\Big(z^{+},z^{-}\Big)}_{L^{4}}^{4}+
\norm{\Lambda^{\alpha}
\Big(|z^+|^2,|z^-|^2\Big)}^2_{L^2}
\]
\begin{eqnarray}\label{L4-estimate}
=-\underbrace{\int_{\mathbb{R}^{3}}\nabla \pi \cdot
(z^{+}\left\vert z^{+}\right\vert ^{2}+z^{-}\left\vert
z^{-}\right\vert
^{2})dx}_{\mathcal{J}_1}+\underbrace{\int_{\mathbb{R}^{3}} (|z^+|^2
z^{+}+|z^-|^2
z^{-}) \cdot (\nabla \times w)dx}_{\mathcal{J}_2} .
\end{eqnarray}%

Taking the operator \mbox{div}, to the first equation of (4.1), and
using the facts $\mbox{div}(\nabla \times w) = 0$, we see
\[
-\Delta \pi = \mbox{div}\mbox{div}(z^+\otimes z^-),
\]
and thus, for $\zeta>1$
\[
||\pi||_{L^{\delta,2}}\leq C||z^+\otimes z^-||_{L^{\delta,2}}\leq
C||z^+||_{L^{2\delta,4}}||z^-||_{L^{2\delta,4}}=C|||z^+|^2||^{1/2}_{L^{\delta,2}}|||z^-|^2||^{1/2}_{L^{\delta,2}}
\]
\[
\leq C(|||z^+|^2||_{L^{\delta,2}}+|||z^-|^2||_{L^{\delta,2}}).
\]
Using integration by parts, H\"{o}lder's inequality in Lorentz space,
Young's inequalities and Sobolev embedding for the fractional power,
we note that for $p>1$, $\delta>2$ and $\frac{5}{2}>\alpha>0$ with
\begin{equation}\label{relation}
\frac{1}{\delta}+\frac{1}{2p}+\frac{5-2\alpha}{6}=1,
\end{equation}
\[
\int_{\mathbb{R}^{3}}\nabla \pi \cdot \left\vert z^{+}\right\vert
^{2}z^{+}\, dx=\int_{\mathbb{R}^{3}}   \nabla (z^{+}\left\vert
z^{+}\right\vert ^{2})\pi\, dx\leq C\int_{\mathbb{R}^{3}}
|z^{+}||\nabla |z^{+}|^2||\pi|\, dx
\]
\[
\leq C\|z^+\|_{L^{2\delta,4}}\|\nabla
|z^+|^2\|_{L^{\frac{6}{5-2\alpha},2}}\||\pi|^{1/2}\|_{L^{2p,\infty}}\||\pi|^{1/2}\|_{L^{2\delta,4}}
\]
\[
=C\||z^+|^2\|^{1/2}_{L^{\delta,2}}\|\nabla
|z^+|^2\|_{L^{\frac{6}{5-2\alpha},2}}\|\pi\|^{1/2}_{L^{p,\infty}}\|\pi\|^{1/2}_{L^{\delta,2}}
\]
\[
\leq C\||z^+|^2\|^{1/2}_{L^{\delta,2}}\|\nabla
|z^+|^2\|_{L^{\frac{6}{5-2\alpha},2}}\|\pi\|^{1/2}_{L^{p,\infty}}\Big(|||z^+|^2||^{1/2}_{L^{\delta,2}}+|||z^-|^2||^{1/2}_{L^{\delta,2}}\Big)
\]
\[
\leq C\||z^+|^2\|_{L^{\delta,2}}\|\nabla
|z^+|^2\|_{L^{\frac{6}{5-2\alpha},2}}\|\pi\|^{1/2}_{L^{p,\infty}}|||z^-|^2||^{1/2}_{L^{\delta,2}}
\]
\[
+C\||z^+|^2\|^{1/2}_{L^{\delta,2}}\|\nabla
|z^+|^2\|_{L^{\frac{6}{5-2\alpha},2}}\|\pi\|^{1/2}_{L^{p,\infty}}|||z^-|^2||^{1/2}_{L^{\delta,2}}
\]
\[
\leq
C\|\pi\|^{1/2}_{L^{p,\infty}}\Big(\||z^+|^2\|_{L^{\delta,2}}+\||z^-|^2\|_{L^{\delta,2}}\Big)\|\nabla
|z^+|^2\|_{L^{\frac{6}{5-2\alpha},2}}.
\]
In the same way, $\int_{\mathbb{R}^{3}}\nabla \pi \cdot
z^{-}\left\vert z^{-}\right\vert ^{2}\, dx$ can be bounded by
\[
C\|\pi\|^{1/2}_{L^{p,\infty}}\Big(\||z^+|^2\|_{L^{\delta,2}}+\||z^-|^2\|_{L^{\delta,2}}\Big)\|\nabla
|z^-|^2\|_{L^{\frac{6}{5-2\alpha},2}}.
\]
And thus, using the Gagliardo–Nirenberg interpolation inequality, it follows
\[
\mathcal{J}_1
\leq
C\|\pi\|^{1/2}_{L^{p,\infty}}\Big(\||z^+|^2\|_{L^{\delta,2}}+\||z^-|^2\|_{L^{\delta,2}}\Big)
\Big(\|\nabla|z^+|^2\|_{L^{\frac{6}{5-2\alpha},2}}+\|\nabla|z^-|^2\|_{L^{\frac{6}{5-2\alpha},2}}\Big)
\]
\[
\leq
C\|\pi\|^{1/2}_{L^{p,\infty}}\Big(\||z^+|^2\|^{1-\Big(\frac{3}{2\alpha}-\frac{3}{q\alpha}\Big)}_{L^{2}}+\||z^-|^2\|^{1-\Big(\frac{3}{2\alpha}-\frac{3}{q\alpha}\Big)}_{L^{2}}\Big)
\]
\[
\times
\Big(\|\nabla|z^+|^2\|^{1+\Big(\frac{3}{2\alpha}-\frac{3}{q\alpha}\Big)}_{L^{\frac{6}{5-2\alpha},2}}+\|\nabla|z^-|^2\|^{1+\Big(\frac{3}{2\alpha}-\frac{3}{q\alpha}\Big)}_{L^{\frac{6}{5-2\alpha},2}}\Big)
\]
\[
\leq C\|\pi\|^{\frac{2\alpha p}{4\alpha
p-2p-3}}_{L^{p,\infty}}\Big(\||z^+|^2\|^2_{L^{2}}+\||z^-|^2\|^2_{L^{2}}\Big)
+\frac{1}{16}\Big(\|\Lambda^{\alpha}|z^+|^2\|^2_{L^2}+\|\Lambda^{\alpha}|z^-|^2\|^2_{L^2}\Big).
\]
For $\mathcal{J}_2$, integrating by parts, we note that
\[
\int_{\mathbb{R}^{3}} |z^+|^2
z^{+} \cdot (\nabla \times w)dx\leq \norm{w}_{L^4}\norm{\nabla
|z^+|^2}_{L^{2}}\norm{z^+}_{L^4}\leq\norm{w}^2_{L^4}\|z^+\|^2_{L^{4}}+\frac{1}{16}\norm{\nabla
|z^+|^2}^2_{L^{2}}
\]
\[
\leq\norm{w}^2_{L^4}\|z^+\|^2_{L^{4}}+\frac{1}{16}\norm{|z^+|^2}^{2\theta}_{L^{2}}\norm{\Lambda^\alpha
|z^+|^2}^{2(1-\theta)}_{L^{2}},\quad \theta=\frac{\alpha-1}{\alpha} 
\]
\[
\leq C(\norm{w}^4_{L^4}+\|z^+\|^4_{L^{4}})+\frac{C}{16}\norm{|z^+|^2}^{2}_{L^{2}}+\frac{1}{16}\norm{\Lambda^\alpha
|z^+|^2}^{2}_{L^{2}}
\]
\[
\leq C(\norm{w}^4_{L^4}+\|z^+\|^4_{L^{4}})+\frac{1}{16}\norm{\Lambda^\alpha
|z^+|^2}^{2}_{L^{2}}.
\]
And thus, $\mathcal{J}_2$ is bounded by
\[
\mathcal{J}_2 \leq  C(\norm{w}^4_{L^4}+\|z^+\|^4_{L^{4}}+\|z^-\|^4_{L^{4}})+\frac{1}{16}\Big(\norm{\Lambda^\alpha
|z^+|^2}^{2}_{L^{2}}+\norm{\Lambda^\alpha
|z^-|^2}^{2}_{L^{2}}\Big).
\]

To get $L^4$-estimate for $w$, as before, multiplying the third equation of \eqref{fmmp-100} by $
\left\vert w\right\vert ^{2}w$, integrating by parts and
summing up, we have
\[
\frac{1}{4}\frac{d}{dt}\norm{w}_{L^{4}}^{4}+
\norm{\Lambda^\gamma
|w|^2}^2_{L^2}+2\chi\norm{w}_{L^{4}}^{4}+\norm{|w|\mbox{div}\
w}_{L^{2}}^{2}
\]
\begin{eqnarray}\label{L4-estimate}
=\underbrace{\frac{\chi}{2}\int_{\mathbb{R}^{3}}|w|^2 w \cdot
(\nabla \times (z^+ + z^-))dx}_{\mathcal{J}_3}-
\underbrace{\int_{\mathbb{R}^{3}} \mbox{div}\ w\ (w \cdot
\nabla|w|^2)dx}_{\mathcal{J}_4}.
\end{eqnarray}%
As same manner as $\mathcal{J}_2$, $\mathcal{J}_3$ is bounded by 
\[
\mathcal{J}_3\leq
C(\norm{w}^{4}_{L^4}+\norm{z^+}^{4}_{L^4}+\norm{z^-}^{4}_{L^4})+\frac{1}{16}\norm{\Lambda^\gamma
|w|^2}^2_{L^2},\quad   \gamma\geq1.
\]
In a similar way, $\mathcal{J}_4$ is also bounded by  
\[
\mathcal{J}_4\leq \kappa\norm{|w|\mbox{div}\ w}_{L^{2}}\norm{\nabla
|w|^2}_{L^{2}}\leq C\norm{\nabla
|w|^2}^2_{L^{2}}+\frac{\kappa}{16}\norm{|w|\mbox{div}\ w}^2_{L^{2}}
\]
\[
\leq
C\norm{w}^4_{L^{4}}+\frac{1}{16}\norm{\Lambda^\gamma|w|^2}_{L^{2}}^2+\frac{\kappa}{16}\norm{|w|\mbox{div}\ w}^2_{L^{2}}, \quad \gamma\geq1.
\]
Let
$Y(t):=\|(z^+,z^-,w)\|^4_{L^4(\mathbb{R}^{3})}$
and thus \eqref{L4-estimate} becomes
\begin{equation}\label{qq111}
 \frac{d}{dt}Y(t)\lesssim (1+ \|\pi\|^{q}_{L^{p,\infty}(\mathbb{R}^{3})})Y(t),\quad\quad q=\frac{2\alpha p}{4\alpha
p-2p-3}.
\end{equation}
Now, we use an argument similar to the one used in the work of Bosia
et al. \cite{BPR14}. For $\kappa > 0$, Choose $q_\kappa = q
-\kappa(q+\frac{\alpha}{2\alpha-1}-\frac{3c_0}{4(2\alpha-1)})$ and
$r_\kappa:=\frac{q-\kappa(q-\frac{\alpha}{2\alpha-1}-\frac{3c_0}{4(2\alpha-1)})}{\frac{2}{3}(q(2\alpha-1)-\alpha)(1-\kappa)+\frac{c_0\epsilon}{2} }$ with
\begin{equation*}\left\{\ba\label{ro}
&\frac{3}{p_{\kappa}}+\frac{2\alpha}{q_{\kappa}}=2(2\alpha-1), \\
&\frac{q_{\kappa}}{p_{\kappa}}=\frac{q\big(1-\kappa\big)}{p}+\frac{c_0\kappa}{2}.
\ea\right.\end{equation*}
Due to the above relation, we get
\[
 \|\pi\|^{q_{\kappa}}_{L^{p_{\kappa},\infty}(\mathbb{R}^{3})}\lesssim
 \|\pi\|^{q(1-\kappa)}_{L^{p ,\infty}(\mathbb{R}^{3})} \|\pi\|^{4\kappa}_{L^{2,\infty}}\lesssim
 \|\pi\|^{q(1-\kappa)}_{L^{p ,\infty}(\mathbb{R}^{3})} \|\pi\|^{4\kappa}_{L^{2 }(\mathbb{R}^{3})}
\]
\begin{equation}\label{2.2000}
\lesssim
 \|\pi\|^{q(1-\kappa)}_{L^{p ,\infty}(\mathbb{R}^{3})} \Big(\norm{|u|^{2}}^{4\kappa}_{L^{2
 }(\mathbb{R}^{3})}+\norm{|b|^{2}}^{4\kappa}_{L^{2
 }(\mathbb{R}^{3})}\Big).
\end{equation}
Since the pair $(p_{\kappa}, q_{\kappa})$ also meets
$3/p_{\kappa}+2\alpha/q_{\kappa}=2(2\alpha-1)$. Using the estimate \eqref{2.2000},
\eqref{qq111} becomes
$$
\frac{d}{dt}Y(t)\lesssim
(1+ \|\pi\|^{q_{\kappa}}_{L^{p_{\kappa},\infty}(\mathbb{R}^{3})})\Big\||u|^{2}\Big\|_{L^{2}(\mathbb{R}^{3})}^{2}
  \lesssim (1+\|\pi\|^{q(1-\kappa)}_{L^{p ,\infty}(\mathbb{R}^{3})} )Y(t)^{1+2\kappa}.
$$
And then integrating with respect to time from $0$ to $t$ with
$0\leq t<T$,
\[
Y(t) \leq CY(0)+C\int_0^t (1+\|\pi\|^{q(1-\kappa)}_{L^{p
,\infty}(\mathbb{R}^{3})} )Y(t)^{1+2\kappa}\,ds,
\]
or equivalently,
\[
\|w^+(t)\|^4_{L^4(\mathbb{R}^{3})}+\|w^-(t)\|^4_{L^4(\mathbb{R}^{3})}
\leq
C\|w_0^+\|^4_{L^4(\mathbb{R}^{3})}+\|w_0^-\|^4_{L^4(\mathbb{R}^{3})}
\]
\[
+C\int_0^t (1+\|\pi\|^{q(1-\kappa)}_{L^{p ,\infty}(\mathbb{R}^{3})})
\|w^+\|^4_{L^4(\mathbb{R}^{3})}+\|w^-\|_{L^4(\mathbb{R}^{3})}^{4(1+2\kappa)}\,ds.
\]
Due to Lemma \ref{grow}, we are now able to complete the proof of
Theorem 1.1 under the assumption $(A)$ in Theorem 1.1.


\it{Part (B):}  Indeed, a proof is almost same to that in the argument in \cite{Duan12} or \cite{Kim22}, however, for the reader's convenience, a sketch of the proof will be given.
Multiplying both side of $\eqref{fmmp-100}_1$ by $z^+|z^+|^{3r-4}$ and then integrating them over $\R^3$
we conclude that
\begin{align*}
\frac{1}{3r-2}\f{d}{dt}&\int_{\mathbb{R}^{3}}|z^+|^{3r-2}dx
+\frac{4(3r-4)}{(3r-2)^2}\int_{\mathbb{R}^{3}}|\Lambda^{\alpha}|
z^+|^{\frac{3r-2}{2}}|^{2}\,dx
\end{align*}
\begin{eqnarray}\label{L4-estimate}
\lesssim \underbrace{\int_{\mathbb{R}^{3}}\nabla \pi \cdot |z^+|^{3r-4}z^+ dx }_{\mathcal{J}_5}+\frac{1}{2}\underbrace{\int_{\mathbb{R}^{3}}|z^+|^{3r-4}z^+ \cdot (\nabla \times w )dx}_{\mathcal{J}_6}
\end{eqnarray}
By the integration by parts and  H\"{o}lder inequality, $\mathcal{J}_5$ is also written by 
\begin{align}\label{aaa-1000}
\mathcal{J}_1&\leq (3r - 4) \int_{\mathbb{R}^{3}} | \pi||\nabla
|z^+||z^+|^{3r-4}\,dx \\
&\leq
\frac{2(3r-4)}{(3r-2)}\Big(\int_{\mathbb{R}^{3}}|\pi|^2|z^+|^{3r-4}\,dx\Big)^\frac{1}{2}
\Big(\int_{\mathbb{R}^{3}}|\nabla|z^+|^{\frac{3r-2}{2}}\,dx\Big)^\frac{1}{2}\notag.
\end{align}
Note that $0\leq I\leq a$ and $0\leq I\leq b$, then $ I\leq
\sqrt{ab}$. Combining  $\mathcal{J}_5$ in \eqref{L4-estimate} and \eqref{aaa-1000}, we get
\begin{align*}
\mathcal{J}_5&\lesssim \Big(\int_{\mathbb{R}^{3}} |\nabla \pi|
|z^+|^{3r-3}\,dx\Big)^{1/2}\Big(\int_{\mathbb{R}^{3}}|\pi|^2|z^+|^{3r-4}\,dx\Big)^\frac{1}{4}
\Big(\int_{\mathbb{R}^{3}}|\nabla|z^+|^{\frac{3r-2}{2}}|\,dx\Big)^\frac{1}{4}\\
&\leq C \Big(\int_{\mathbb{R}^{3}} \Big(|\nabla \pi|\Big(
|z^+|^{2}+|z^-|^{2}\Big)^{(3r-3)/2}\,dx\Big)^{2/3}\Big(\int_{\mathbb{R}^{3}}\Big(|\pi|\Big(
|z^+|^{2}+|z^-|^{2}\Big)^{(3r-4)/2}\,dx\Big)^{1/3}\\
&\hspace{2cm}+\frac{3r-4}{(3r-2)^2}\Big(\int_{\mathbb{R}^{3}}|\nabla
|z^+|^{\frac{3r-2}{2}}|^2\,dx\Big).
\end{align*} Due to
\[
\int_{\mathbb{R}^{3}} | \pi|^2\Big(
|z^+|^{2}+|z^-|^2\Big)^{\frac{3r-4}{2}}\,dx\lesssim
\|\pi\|_{L^{\frac{3r}{2},6r-6}}^2\||z^+|^2+|z^-|^2\|_{L^{\frac{3r}{2},\frac{3r-3}{2}}}^{\frac{3r-4}{2}}
\]
\[
\lesssim
\||z^+|^2+|z^-|^2\|_{L^{\frac{3r}{2},6r-6}}^2\||z^+|^2+|z^-|^2\|_{L^{\frac{3r}{2},\frac{3r-3}{2}}}^{\frac{3r-4}{2}}
\]
\[
\lesssim
\||z^+|^2+|z^-|^2\|_{L^{\frac{3r}{2},\frac{3r-3}{2}}}^2\||z^+|^2+|z^-|^2\|_{L^{\frac{3r}{2},\frac{3r-3}{2}}}^{\frac{3r-4}{2}}=\||z^+|^2+|z^-|^2\|_{L^{\frac{3r}{2},\frac{3r-3}{2}}}^{\frac{3r}{2}},
\]
and
\[
\int_{\mathbb{R}^{3}} | \nabla \pi|\Big(
|z^+|^{2}+|z^-|^2\Big)^{\frac{3r-3}{2}}\,dx\lesssim \|\nabla
\pi\|_{L^{r,\infty}}\||z^+|^2+|z^-|^2\|_{L^{\frac{3r}{2},\frac{3r-3}{2}}}^{\frac{3r-3}{2}},
\]
$\mathcal{J}_5$ is estimated by
\begin{equation}\label{aaaaa-10}
\mathcal{J}_5\leq C\|\nabla
\pi\|^{\frac{2}{3}}_{L^{r,\infty}}\||z^+|^2+|z^-|^2\|_{L^{\frac{3r}{2},\frac{3r-2}{2}}}^{\frac{3r-2}{2}}+\frac{3r-4}{(3r-2)^2}\Big(\int_{\mathbb{R}^{3}}|\nabla
|z^+|^{\frac{3r-2}{2}}|^2\,dx\Big).
\end{equation}
Next, for $\mathcal{J}_6$, using the H\"{o}lder and Young inequalities, we have
\begin{equation}\label{aaaaa-20}
\mathcal{J}_6\lesssim \|w\|_{L^{3r-2}}\||z^+|^{\frac{3r-4}{2}}\|_{L^{\frac{2(3r-2)}{3r-4}}}\|\nabla |z^+|^{\frac{3r-2}{2}}\|_{L^2}
\end{equation}
\[
\lesssim\|w\|^2_{L^{3r-2}}\||z^+|^{\frac{3r-4}{2}}\|^2_{L^{\frac{2(3r-2)}{3r-4}}}+\|\nabla |z^+|^{\frac{3r-2}{2}}\|^2_{L^2}
\lesssim(\|w\|^{3r-2}_{L^{3r-2}}+\||z^+\|^{3r-2}_{L^{3r-2}})+\|\nabla |z^+|^{\frac{3r-2}{2}}\|^2_{L^2}.
\]
And them, considering the estimates \eqref{aaaaa-10} and \eqref{aaaaa-20}, \eqref{L4-estimate} reduces
\begin{equation}
\begin{aligned}\label{2.13}
\f{d}{dt}&\int_{\mathbb{R}^{3}}|z^+|^{3r-2}dx
+\int_{\mathbb{R}^{3}}|\Lambda^{\alpha}|
z^+|^{\frac{3r-2}{2}}|^{2}\,dx
\end{aligned}
\end{equation}
\[
\lesssim \|\nabla
\pi\|^{\frac{2}{3}}_{L^{r,\infty}}\||z^+|^2+|z^-|^2\|_{L^{\frac{3r}{2},\frac{3r-2}{2}}}^{\frac{3r-2}{2}}+ C(\|w\|^{3r-2}_{L^{3r-2}}+\||z^+\|^{3r-2}_{L^{3r-2}})+\|\nabla |z^+|^{\frac{3r-2}{2}}\|^2_{L^2}.
\]
\[
\leq C \|\nabla \pi\|^{\frac{2}{3}}_{L^{r,\infty}}\||z^+|^2+|z^-|^2\|_{L^{\frac{3r}{2},\frac{3r-2}{2}}}^{\frac{3r-2}{2}}+ C(\|w\|^{3r-2}_{L^{3r-2}}+\||z^+\|^{3r-2}_{L^{3r-2}})+\frac{1}{256}\|\Lambda^{\alpha} |z^+|^{\frac{3r-2}{2}}\|^2_{L^2}.
\]
where we use the estimate
\[
\|\nabla |z^+|^{\frac{3r-2}{2}}\|^2_{L^2}
\leq C \|z^+\|^{3r-2}_{L^{3r-2}}\|\Lambda^{\alpha} |z^+|^{\frac{3r-2}{2}}\|^{2(1-\theta)}_{L^2}
\leq  \||z^+|^{\frac{3r-2}{2}}\|^{2\theta}_{L^2}+\frac{1}{256}\|\Lambda^{\alpha} |z^+|^{\frac{3r-2}{2}}\|^{2}_{L^2}.
\]
 In a similar fashion, if you do it for the equation
$\eqref{fmmp-100}_2$, we have
\begin{equation}
\begin{aligned}\label{2.14}
\f{d}{dt}&\int_{\mathbb{R}^{3}}|z^-|^{3r-2}dx
+\int_{\mathbb{R}^{3}}|\Lambda^{\alpha}|
z^-|^{\frac{3r-2}{2}}|^{2}\,dx
\end{aligned}
\end{equation}
\[
\leq C \|\nabla \pi\|^{\frac{2}{3}}_{L^{r,\infty}}\||z^+|^2+|z^-|^2\|_{L^{\frac{3r}{2},\frac{3r-2}{2}}}^{\frac{3r-2}{2}}+ C(\|w\|^{3r-2}_{L^{3r-2}}+\||z^-\|^{3r-2}_{L^{3r-2}})+\frac{1}{256}\|\Lambda^{\alpha} |z^-|^{\frac{3r-2}{2}}\|^2_{L^2}.
\]
After summing up \eqref{2.13} and \eqref{2.14}, using Sobolev
embedding and Young's inequality, we obtain
\begin{equation}
\begin{aligned}\label{2.15-10}
\frac{d}{dt}&\int_{\mathbb{R}^{3}}\Big(|z^+|^{3r-2}+|z^-|^{3r-2}\Big)dx
+\int_{\mathbb{R}^{3}}\Big(|\Lambda^{\alpha}|
z^+|^{\frac{3r-2}{2}}|^{2}+|\Lambda^{\alpha}|
z^-|^{\frac{3r-2}{2}}|^{2}\Big)\,dx\\
&\lesssim\|\nabla
\pi\|^{\frac{2}{3}}_{L^{r,\infty}}\Big(\||z^+|^{\frac{3r-2}{2}}\|_{L^{\frac{6r}{3r-2}},1}^{2}+\||z^-|^{\frac{3r-2}{2}}\|_{L^{\frac{6r}{3r-2}},1}^{2}\Big)\\
&\hspace{4cm}+C(\|w\|^{3r-2}_{L^{3r-2}}+\||z^+\|^{3r-2}_{L^{3r-2}}+\||z^-\|^{3r-2}_{L^{3r-2}})\\
&\lesssim\|\nabla
\pi\|^{\frac{2}{3}}_{L^{r,\infty}}\Big(\||z^+|^{\frac{3r-2}{2}}\|_{L^{2}}^{(2-\frac{2}{r\alpha})}\|\Lambda^{\alpha}|z^+|^{\frac{3r-2}{2}}\|_{L^{2}}^{\frac{2}{r\alpha}}
+\||z^-|^{\frac{3r-2}{2}}\|_{L^{2}}^{(2-\frac{2}{r\alpha})}\|\Lambda^{\alpha}|z^-|^{\frac{3r-2}{2}}\|_{L^{2}}^{\frac{2}{r\alpha}}\Big)\\
&\hspace{4cm}+C(\|w\|^{3r-2}_{L^{3r-2}}+\||z^+\|^{3r-2}_{L^{3r-2}}+\||z^-\|^{3r-2}_{L^{3r-2}})\\
&\lesssim\|\nabla
\pi\|^{\frac{2r\alpha}{3r\alpha-3}}_{L^{r,\infty}}\Big(\|z^+|^{\frac{3r-2}{2}}+|z^-|^{\frac{3r-2}{2}}\|^2_{L^2}\Big)
+\frac{1}{8}\int_{\mathbb{R}_+^{3}}\Big(|\Lambda^{\alpha}|
z^+|^{\frac{3r-2}{2}}|^{2}+|\Lambda^{\alpha}|
z^-|^{\frac{3r-2}{2}}|^{2}\Big)\,dx\\
&\hspace{4cm}+C(\|w\|^{3r-2}_{L^{3r-2}}+\|z^+\|^{3r-2}_{L^{3r-2}}+\|z^-\|^{3r-2}_{L^{3r-2}})\\
&\lesssim\|\nabla
\pi\|^{\frac{2r\alpha}{3r\alpha-3}}_{L^{r,\infty}}\Big(\|z^+\|^{3r-2}_{L^{3r-2}(\mathbb{R}^{3})}+\|z^-\|^{3r-2}_{L^{3r-2}}\Big)
+\frac{1}{8}\int_{\mathbb{R}^{3}}\Big(|\Lambda^{\alpha}|
z^+|^{\frac{3r-2}{2}}|^{2}+|\Lambda^{\alpha}|
z^-|^{\frac{3r-2}{2}}|^{2}\Big)\,dx\\
&\hspace{4cm}+C(\|w\|^{3r-2}_{L^{3r-2}}+\|z^+\|^{3r-2}_{L^{3r-2}}+\|z^-\|^{3r-2}_{L^{3r-2}}).
\end{aligned}
\end{equation} 
To control $L^{3r-2}$-estimate for $w$, multiplying both side of $\eqref{fmmp-100}_3$ by  $w|w|^{3r-4}$, we have 
\begin{align*}
\frac{1}{3r-2}\f{d}{dt}&\int_{\mathbb{R}^{3}}|w|^{3r-2}dx
+\frac{4(3r-4)}{(3r-2)^2}\int_{\mathbb{R}^{3}}|\Lambda^{\gamma}|
w|^{\frac{3r-2}{2}}|^{2}\,dx
\end{align*}
\[
+\int_{\mathbb{R}^{3}}|w|^{3r-2}dx+\int_{\mathbb{R}^{3}}|w|^{3r-4}|\mbox{\rm div}\ w|^2\,dx
\]
\begin{eqnarray}\label{2.16}
=\underbrace{\frac{\chi}{2}\int_{\mathbb{R}^{3}}|w|^{3r-4}w \cdot
(\nabla \times (z^++z^-) )dx}_{\mathcal{J}_7}+
\underbrace{\int_{\mathbb{R}^{3}} \mbox{\rm div}\ w \cdot w\  \mbox{\rm div}(|w|^{3r-4})dx}_{\mathcal{J}_8}.
\end{eqnarray}
As same manner as $\mathcal{J}_2$, the term $\mathcal{J}_7$ is bounded by
\begin{equation}\label{2.17}
\mathcal{J}_7\lesssim(\|w\|^{3r-2}_{L^{3r-2}}+\||z^+\|^{3r-2}_{L^{3r-2}}+\||z^-\|^{3r-2}_{L^{3r-2}})+\frac{1}{256}\|\Lambda^{\gamma} |w|^{\frac{3r-2}{2}}\|^2_{L^2}.
\end{equation}
For $\mathcal{J}_8$, we get
\begin{equation}\label{2.18}
\mathcal{J}_8\leq C\|\nabla |w|^{\frac{3r-2}{2}}\|^2_{L^2}+\frac{1}{256}\int_{\mathbb{R}^{3}}|w|^{3r-4}|\mbox{\rm div}\ w|^2\,dx
\end{equation}
\[
\leq C\|w\|^{3r-2}_{L^{3r-2}}+\frac{1}{256}\Big(\|\Lambda^{\gamma} |w|^{\frac{3r-2}{2}}\|^2_{L^2}+\||w|^{\frac{3r-4}{2}}|\mbox{\rm div}\ w|\|^2_{L^2}\Big).
\]
Summing up all estimates \eqref{2.15-10}--\eqref{2.18}, we have
\begin{equation*}
\begin{aligned}\label{2.15}
\frac{d}{dt}&(\|z^+\|^{3r-2}_{L^{3r-2}(\mathbb{R}^{3})}+\|z^-\|^{3r-2}_{L^{3r-2}(\mathbb{R}^{3})}+\|w\|^{3r-2}_{L^{3r-2}(\mathbb{R}^{3})})\\
&\lesssim\|\nabla
\pi\|^{\frac{2r\alpha}{3r\alpha-3}}_{L^{r,\infty}}\Big(\|z^+\|^{3r-2}_{L^{3r-2}(\mathbb{R}^{3})}+\|z^-\|^{3r-2}_{L^{3r-2}}\Big)
+(\|w\|^{3r-2}_{L^{3r-2}}+\|z^+\|^{3r-2}_{L^{3r-2}}+\|z^-\|^{3r-2}_{L^{3r-2}}).
\end{aligned}
\end{equation*} 
Let
$\mathcal{Y}(t):=\|z^+\|^{3r-2}_{L^{3r-2}(\mathbb{R}^{3})}+\|z^-\|^{3r-2}_{L^{3r-2}(\mathbb{R}^{3})}+\|w\|^{3r-2}_{L^{3r-2}(\mathbb{R}^{3})}$
and then \eqref{2.15} becomes
\[
\mathcal{Y}(t) \leq C\| \nabla \pi
\|^{\f{2r\alpha}{3r\alpha-3}}_{L^{r,\infty}(\mathbb{R}^{3})}\mathcal{Y}(t)+\mathcal{Y}(t).
\]
As the previous way, it allow us to finish the proof of Theorem
\ref{the1.1}.
\end{pf}

\section{Proof of Theorem \protect\ref{the1.2}}
For this, according to the argument in \cite{Zhou06} or \cite{JWW20}, we can establish a Serrin's type regularity criterion on the gradient of pressure  function $\pi$. Indeed, from
\eqref{L4-estimate}, we know
\[
\frac{1}{4}\frac{d}{dt}\|(u,b,w)\|_{L^{4}}^{4}+
\|\nabla(|u|^2,|b|^2,|w|^2)\|^2_{L^2}
\]
\[
+\norm{|u||\nabla u|}^2_{L^2}+\norm{|b||\nabla b|}^2_{L^2}+\norm{|w||\nabla w|}^2_{L^2}+2\chi\|w\|_{L^{4}}^{4}+|||w|\mbox{div}\
w||_{L^{2}}^{2}
\]
\[
\lesssim \underbrace{\int_{\mathbb{R}^{3}}|\nabla \mathcal{P}| |u|^3|\nabla |u|^2| dx}_{\mathcal{J}_1}
+\underbrace{\int_{\mathbb{R}^{3}}(b\cdot \nabla) b\cdot |u|^2u dx}_{\mathcal{J}_2}
\]
\[
+\frac{1}{2}\underbrace{\int_{\mathbb{R}^{3}}\nabla(|b|^2)\cdot |u|^2u \,dx}_{\mathcal{J}_3}+\underbrace{\int_{\mathbb{R}^{3}}(b\cdot \nabla)u\cdot |b|^2b dx dt}_{\mathcal{J}_4}+\frac{\chi}{2}\underbrace{\int_{\mathbb{R}^{3}}|w|^2 w \cdot (\nabla \times u )dx}_{\mathcal{J}_5}
\]
\begin{eqnarray}\label{L4-estimate-10}
+\underbrace{\frac{\chi}{2}\int_{\mathbb{R}^{3}}|w|^2 w \cdot
(\nabla \times u )dx}_{\mathcal{J}_6}-
\underbrace{\int_{\mathbb{R}^{3}} \mbox{div}\ w\ (w \cdot
\nabla|w|^2)dx}_{\mathcal{J}_7}
\end{eqnarray}%
For the result,  $\mathcal{J}_1$ only has been changed as follows: for $r>1$
\[
\int_{\mathbb{R}^{3}}\nabla \mathcal{P} \cdot|u|^{2}u\, dx\leq \||\nabla \mathcal{P}|^{1/2}\|_{L^{4,4}} \||\nabla
\mathcal{P}|^{1/2}\|_{L^{2r,\infty}}
\||u^3\|_{L^{\frac{4r}{3r-2},\frac{4}{3}}}
\]
\[
=\|\nabla\mathcal{P}\|^{1/2}_{L^{2,2}} \|\nabla
\mathcal{P}\|^{1/2}_{L^{r,\infty}}
\|u\|^3_{L^{\frac{12r}{3r-2},4}}\leq
\frac{1}{4}\|\nabla \pi\|^2_{L^2}+ C\|\nabla
\mathcal{P}\|^{3/2}_{L^{r,\infty}}
\|u\|^4_{L^{\frac{12r}{3r-2},4}}
\]
\[
\leq \frac{1}{4}\|\nabla \mathcal{P}\|^2_{L^2}+ C\|\nabla
\mathcal{P}\|^{\frac{2}{3}}_{L^{\gamma,\infty}}
\||u|^2\|^2_{L^{\frac{6r}{3r-2},2}}
\]
\[
\leq \frac{1}{4}\|\nabla \mathcal{P}\|^2_{L^2}+ C\|\nabla
\mathcal{P}\|^{\frac{2}{3}}_{L^{r,\infty}}
\||u|^2\|^{2(1-\frac{1}{r})}_{L^{2,2}}\norm{\nabla|u|^2}^{\frac{2}{r}}_{L^{2,2}}
\]
\[
\leq \frac{1}{16}\|\nabla
\mathcal{P}\|^2_{L^2}+\frac{1}{8}\norm{\nabla|u|^2}^{2}_{L^{2}}
+C\|\nabla \mathcal{P}\|^{\frac{2r}{3(r-1)}}_{L^{r,\infty}}
\|u\|^{4}_{L^{4}}, \quad
\]
and thus
\[
\mathcal{J}_1 \leq \frac{1}{4}\|\nabla
\mathcal{P}\|^2_{L^2}+\frac{1}{16}\norm{\nabla|u|^2}^{2}_{L^{2}}
+C\|\nabla \mathcal{P}\|^{\frac{2r}{3(r-1)}}_{L^{r,\infty}}
\norm{u}^{4}_{L^{4}}.
\]
Using the following estimate,
\[
\|\nabla v\|^2_{L^2}\lesssim \|(u\cdot \nabla
)u+(b\cdot \nabla )b\|^2_{L^2}\lesssim  \norm{|u||\nabla u|}^2_{L^2}+ \norm{|b||\nabla b|}^2_{L^2}.
\]
we get
\[
\mathcal{J}_1 \leq C \|\nabla
\mathcal{P}\|^{\frac{2r}{3(r-1)}}_{L^{r,\infty}}\norm{u}^{4}_{L^{4}}+\frac{1}{8}( \norm{|u||\nabla u|}^2_{L^2}+ \norm{|b||\nabla b|}^2_{L^2}).
\]
Using the integration by parts,  $\mathcal{J}_2$,  $\mathcal{J}_3$ and $\mathcal{J}_4$ is bounded by
\[
\int_0^T\int_{\mathbb{R}^{3}}|u||b|^2 (\nabla |u|^2+\nabla |b|^2) dx dt
\leq C(\norm{(|u|^2+|b|^2)b}^2_{L^2}+\frac{1}{16}(\norm{\nabla |u|^2}_{L^2}^2+\norm{\nabla |b|^2}_{L^2}^2)
\]
\[
\leq  C\norm{b}^{\frac{2a_1}{a_1-3}}_{L^{a_1}}(\norm{ |u|^2}^{2}_{L^2}+\norm{ |b|^2}^{2}_{L^2})+\frac{1}{16}(\norm{\nabla|u|^2}^{2}_{L^2}+\norm{\nabla|b|^2}^{2}_{L^2})
\]
where we use the following inequality:
\[
\norm{b}^2_{L^{a_1}}\norm{|u|^2}^2_{L^{\frac{2a_1}{a_1-2}}}\lesssim  \norm{b}^2_{L^{a_1}}\norm{|u|^2}^{2(1-\frac{3}{a_1})}_{L^2}\norm{\nabla|u|^2}^{\frac{6}{a_1})}_{L^2}
\leq C \norm{b}^{\frac{2a_1}{a_1-3}}_{L^{a_1}}\norm{ |u|^2}^{2}_{L^2}+\frac{1}{16}\norm{\nabla|u|^2}^{2}_{L^2}
\]
In a similar way, for $\mathcal{J}_5$ and $\mathcal{J}_6$, it shows
\[
|J_5|\leq C\int_{\mathbb{R}^{3}}(|u|^{4}+|w|^{4})\,dx+\frac{1}{16}\int_{\mathbb{R}^{3}}||w||\nabla w||^2\,dx,
\]
and
\[
|J_6|\leq C\int_{\mathbb{R}^{3}} |\mbox{div}\ w|^2dx+\frac{1}{16}\int_{\mathbb{R}^{3}} |\nabla |w|^2|^2dx.
\]
Plugging this into \eqref{L4-estimate-10}, we get
$$\ba
\frac{d}{dt}&\|(u,b,w)\|_{L^{4}}^{4}+
\|\nabla(|u|^2,|b|^2,|w|^2)\|^2_{L^2}
+\norm{|u||\nabla u|}^2_{L^2}+\norm{|b||\nabla b|}^2_{L^2}+\norm{|w||\nabla w|}^2_{L^2}\\
&\lesssim \| \nabla  \mathcal{P} \|^{p}_{L^{r,\infty}(\mathbb{R}^{3})}\|(u,b,w)\|_{L^{4}}^{4}
+\norm{b}^{\frac{2a_1}{a_1-3}}_{L^{a_1}}\|(u,b,w)\|_{L^{4}}^{4},\quad q=\frac{2r}{3(r-1)}\\
 &\lesssim\| \nabla  \mathcal{P} \|^{p_{\kappa}}_{L^{r_{\kappa},\infty}(\mathbb{R}^{3})}\|(u,b,w)\|_{L^{4}}^{4}+\norm{b}^{\frac{2a_1}{a_1-3}}_{L^{a_1}}\|(u,b,w)\|_{L^{4}}^{4}\\
 &\lesssim\| \nabla \mathcal{P} \|^{p(1-\kappa)}_{L^{r,\infty }(\mathbb{R}^{3})}\| \nabla \Pi \|^{c_{1}\kappa}_{L^{2 }(\mathbb{R}^{3})}\|(u,b,w)\|_{L^{4}}^{4}+\norm{b}^{\frac{2a_1}{a_1-3}}_{L^{a_1}}\|(u,b,w)\|_{L^{4}}^{4}\\
 &\lesssim\| \nabla  \mathcal{P} \|^{p(1-\kappa)}_{L^{r,\infty}(\mathbb{R}^{3})}\Big\| |u||\nabla u|+|b||\nabla b| \Big\|^{c_{1}\kappa}_{L^{2 }(\mathbb{R}^{3})}\|(u,b,w)\|_{L^{4}}^{4}\\ &\leq C\| \nabla \mathcal{P} \|^{\f{2p(1-\kappa)}{2-c_{1}\kappa}}_{L^{q,\infty}(\mathbb{R}^{3})} \| u \|^{\f{8}{2-c_{1}\kappa}}_{L^{4 }(\mathbb{R}^{3})}+\norm{b}^{\frac{2a_1}{a_1-3}}_{L^{a_1}}\|(u,b,w)\|_{L^{4}}^{4}\\
&\quad \quad\quad \quad\quad \quad +\f18\Big(\norm{ |u||\nabla u| }^{2}_{L^{2 }(\mathbb{R}^{3})}+\norm{|b||\nabla b|}^{2}_{L^{2 }(\mathbb{R}^{3})}\Big)\\
&\leq C\| \nabla \mathcal{P} \|^{p(1-\delta)}_{L^{r,\infty}(\mathbb{R}^{3})} \|(u,b,w)\|_{L^{4}(\mathbb{R}^{3})}^{4(1+2\delta) }+\norm{b}^{\frac{2a_1}{a_1-3}}_{L^{a_1}}\|(u,b,w)\|_{L^{4}}^{4}\\
&\quad \quad\quad \quad\quad \quad +\f18\Big(\norm{ |u||\nabla u| }^{2}_{L^{2 }(\mathbb{R}^{3})}+\norm{|b||\nabla b|}^{2}_{L^{2 }(\mathbb{R}^{3})}\Big)\ea$$
Notice that $2/p_{\kappa}+3/r_{\kappa}=3$. Chooing $\delta=\f{(2-c_{1})\kappa}{2-c_{1}\kappa},~~c_{1}=\f{4}{3}$, it finally follows
 $$
\frac{d}{dt}\|(u,b,w)\|_{L^{4}}^{4}\lesssim \| \nabla \mathcal{P} \|^{p(1-\delta)}_{L^{q,\infty}(\mathbb{R}^{3})}  \|(u,b,w)\|_{L^{4}(\mathbb{R}^{3})}^{4(1+2\delta) }+\norm{b}^{\frac{2a_1}{a_1-3}}_{L^{a_1}}\|(u,b,w)\|_{L^{4}}^{4}.$$
As the previous way, it allow us to finish the proof of Theorem
\ref{the1.2}.

 \section*{Acknowledgments}
We would like to appreciate the anonymous referee for valuable
comments.
Jae-Myoung Kim was
supported by National Research Foundation of Korea Grant funded by
the Korean Government (NRF-2020R1C1C1A01006521).

\end{document}